\documentclass[12pt]{amsart}
\usepackage{amssymb}
\usepackage{enumerate}

\newtheorem{fact}{Fact}[section]
\newtheorem{lemma}[fact]{Lemma}
\newtheorem{theorem}[fact]{Theorem}

\newtheorem{rremark}[fact]{Remark}
\newtheorem{exa}[fact]{Example}
 \newtheorem{prob}[fact]{Problem}
 \newtheorem{disc}[fact]{Discussion}
\newtheorem{deta}[fact]{Details}
\newtheorem{remarks}[fact]{Remarks}
\newtheorem{corollary}[fact]{Corollary}

\newtheorem{proposition}[fact]{Proposition}
\newenvironment{remark}{\begin{rremark} \rm}{\end{rremark}}
\newenvironment{example}{\begin{exa} \rm}{\end{exa}}
\newenvironment{discussion}{\begin{disc} \rm}{\end{disc}}

    \def\sqr#1#2{{\vcenter{\hrule height.#2pt
        \hbox{\vrule width.#2pt height#1pt \kern#1pt
            \vrule width.#2pt}\hrule height.#2pt}}}
    \def\square{\mathchoice\sqr67\sqr67\sqr{2.1}6\sqr{1.5}6}

\def\qed{~\hfill$\square$}

\newcommand{\dd}{\partial}
\newcommand{\C}{{\mathbb C}}
\newcommand{\Z}{{\mathbb Z}}
\newcommand{\Q}{{\mathbb Q}}

\newcommand{\p}{{\mathbb P}}
\DeclareMathOperator{\tp}{Tp}
\DeclareMathOperator{\ptp}{\p Tp}

\DeclareMathOperator{\codim}{codim}

\DeclareMathOperator{\Ker}{Ker}

\begin{document}

\begin{abstract}

 Coincident root loci are subvarieties of $S^d\C^2$ --- the space of binary forms of degree $d$ ---
labelled by partitions of $d$. Given a partition $\lambda$,  let
$X_\lambda$ be the set of forms with root multiplicity
corresponding to $\lambda$. There is a natural action of
$GL_2(\C)$ on $S^d\C^2$ and the coincident root loci are invariant
under this action. We calculate their equivariant Poincar\'e duals
generalizing formulas of Hilbert and Kirwan. In the second part we
apply these results to present the cohomology ring of the
corresponding moduli spaces (in the GIT sense) by geometrically
defined relations.

\end{abstract}

\title{Coincident root loci of binary forms}
\author{L. M. Feh\'er}
\address{Department of Analysis, Eotvos University, Budapest}
\email{lfeher@math-inst.hu}
\author{A. N\'emethi}
\address{Department of Mathematics, The Ohio State University}
\email{nemethi@math.ohio-state.edu}
\author{R. Rim\'anyi}
\address{Department of Mathematics, University of North Carolina at Chapel Hill}
\email{rimanyi@email.unc.edu}

\thanks{\noindent Supported by FKFP0055/2001 (1st and 3rd author),
        OTKA T029759 (3rd author) and NSF grant DMS-0088950 (2nd author)
AMS Subject classification 14N10, 57R45}
\maketitle

\section{Introduction}\label{intr}

One of the main goals of Geometric Invariant Theory is to
calculate the cohomology ring of a geometric quotient. In the case
when all semistable point are stable several techniques were
developed. But even for very simple representations this condition
is not satisfied. In this paper we study the action of $GL(2)$ on
the space of binary forms in degree $d$. In the odd case methods
of \cite{kirwan}, \cite{jeffrey-kirwan:nonabelian}, \cite{martin}
can be applied, but none of these methods work in the even case.
We show how equivariant Poincar\'e-dual calculations lead to
relations for the cohomology ring in both the odd and the even
case. These equivariant Poincar\'e-dual (a.c.a. Thom polynomial)
calculations are also interesting on their own right since they
generalize formulas of Hilbert and Kirwan on coincident root loci.
These calculations don't only lead to explicit relations for these
cohomology rings but also identify them with the equivariant
Poincar\'e-duals of the simplest unstable coincident root loci.

\medskip

Consider the $d$-th symmetric power $S^d\C^2$ of the standard representation of $GL_2(\C)$, that is the action of $G$ on the space $V_d$ of degree $d$ homogeneous polynomials in two variables $x$, $y$. For any  partition $\lambda=(\lambda_1,\lambda_2,\ldots,\lambda_n)$ of $d$ (i.e. $\sum_j \lambda_j=d$) we define
$$X_\lambda=\{B(x,y)\in V_d\ |\ B=\prod_{j=1}^n L_j^{\lambda_j}\ \text{for some linear forms}\ L_j\},$$
which is a subvariety  invariant under the group action.
It is called the {\em coincident root loci} associated with $\lambda$.
Clearly,  it is a cone in $V_d$, let $\p X_\lambda$  be its projectivization in the
projective space $\p V_d$.   In fact, it is more convenient to
use a different notation for partitions: $\lambda=(1^{e_1}2^{e_2}\ldots r^{e_r})$ will mean the partition consisting of $e_1$ copies of 1, $e_2$ copies of 2, etc.  Then  $\sum ie_i=d$,   $\sum e_i=n$ and
 the complex dimension of  $\p X_\lambda$ is exactly $n$.

The study of coincident root loci probably started with Cayley.  E.g.,
the  very first question of this type asks the characterization of
polynomials $B$  with a double root. The answer is clearly the vanishing of the
discriminant  which provides in this way an equation for
$X_{(1^{d-2}2)}$. For higher codimensional coincident root loci finding the defining equations is very complicated (see \cite{chipalkatti} for recent results). However, important   geometric
information can be obtained  about these subvarieties.  E.g.,
 the starting point of the present
 paper was Hilbert's formula which calculates the degree of
$\p X_\lambda\subset \p V_d$:
$$\text{deg}(\p X_\lambda)=\frac{n!}{\prod_i (e_i!)}\prod_i i^{e_i}.$$
We can interpret this formula as follows:  for a generic family  of polynomials parametrized by  a
projective space of dimension equal to the codimension of $\p X_\lambda$  the number of
 polynomials in the family with root multiplicity $\lambda$ is $\text{deg}(\p X_\lambda)$.

Generalizing this we arrive to the theory of degeneracy loci. Suppose we have a vector bundle $E\to M$
with fiber $S^d\C^2$ and a generic section $s:M\to E$. Let $s^{-1}(X_\lambda)$ be the set of points in
$M$ where the value of $s$ is in $X_\lambda$.  Its Poincar\'e dual $[s^{-1}(X_\lambda)]\in H^*(M)$
measures  the ``size" of $s^{-1}(X_\lambda)$. It turns out that for any $S^d\C^2$-bundle,
$[s^{-1}(X_\lambda)]$ can be deduced  from the corresponding cohomology class
 of the universal  bundle
associated with the $GL_2(\C)$--representation  $S^d\C^2$. This universal invariant is called
the {\em $GL_2(\C)$-equivariant Poincar\'e dual},
or {\em  Thom polynomial }  of $X_\lambda$ in $S^d\C^2$.
In section 3 we determine all these polynomials.

Calculating equivariant Poincar\'e duals for invariant subvarieties of representations has a long history. We can interpret many results of the nineteenth century algebraic geometers in these terms. From the 1970's the main method was a type of resolution of the subvariety,
 initiated by Porteous \cite{port}. The method requires a deep understanding of the geometry of the resolution and can be carried out only in special cases. Most examples can be found in \cite{fulton}. The first and third author designed a different method (the method of restriction equations, see \cite{cr}) based on ideas coming from calculating Thom polynomials in singularity theory \cite{thompol}. However the method of restriction equations works well mainly if the representation has finitely many orbits which is usually not the case (e.g. for $S^d\C^2$ if $d>3$).

In this paper we return to the technique of resolution, however in a very different way.  The main novelty
is that our new approach  requires only  knowledge of some basic  cohomological data.
Consequently, the method is more flexible. We illustrate  this method here by the coincident root loci,
but the range of applications is much wider. (For example, in a forthcoming paper we plan to discuss
the case  of loci of reducible hypersurfaces.)

Parallel to our work B. K\H om\H uves also provided a presentation of  these Poincar\'e duals in
a completely  different form \cite{balazs}. He worked more in the spirit of the method of restriction
equations studying incidences of the coincident root loci with the orbits $X_{(i,d-i)}$.
We are in the process of comparing these formulas.

In section \ref{moduli} we study the cohomology ring of the moduli
space of the representation $S^d\C^2$ (in the Geometric Invariant
Theory sense). Following the paper of Atiyah and Bott \cite{ab} a
whole theory for calculating cohomology rings of the moduli space
of representations was built up by F. Kirwan; as well as more
algebraic methods were successfully applied by e.g. M. Brion
\cite{brion}, S. Martin \cite{martin}. However, the application of
the general theorems to specific examples is often not easy. Our
approach results explicit presentations of the rational cohomology
rings $H^*_G(X^{ss})$, $H^*(X^{ss}/G)$ and $H^*_G(X^{s})\cong
H^*(X^{s}/G)$ in terms of generators and relations (if $d$ is odd
then all these rings coincide, but for the even case they are
different). We wish to emphasize that a main advantage of our
presentation of the cohomology rings is that we attribute to the
set of relations deep geometric significance: they are the
universal Thom polynomials of some distinguished spaces
$X_\lambda$.

\section{Review on affine and projective Thom polynomials}\label{tp}

Let the group $G$ act on the complex vector space $V$,
and let $\eta$ be an invariant variety in $V$, which supports
a fundamental class (for more details  see \cite{cr}).
Then define the (affine) {\em Thom polynomial}  of $\eta$  as the
Poincar\'e dual of the fundamental homology class
of $\eta$ in equivariant cohomology:
$$\tp_\eta=\text{Poincar\'e dual of }\ [\eta]\in H^*_G(V,\Z).$$
The vector space $V$ is contractible,  hence  the ring $H_G^*(V,\Z)$ is naturally
isomorphic to $H^*(BG,\Z)$, the ring
of $G$-characteristic classes. The degree of $\tp_\eta$
is the real codimension $2c$ of $\eta$ in $V$, hence  $\tp_\eta\in H^{2c}(BG,\Z)$.
The direct geometric meaning of $\tp_\eta$ is the following.

Consider a
fiber bundle $\xi$ with fiber $V$ and structure group $G$ over a manifold $M$.
Because of its invariance, the set $\eta$ can be defined in each fiber, let
the union of these be $\eta(\xi)$. Then consider those points where a generic
section $s$ of $\xi$ hits $\eta$, that is $s^{-1}(\eta(\xi))\subset M$.
By Poincar\'e duality this set defines a cohomology class in $M$.
Standard arguments show that this class  equals $\tp_\eta(\xi):=
f_\xi^*\tp_\eta$, where $f_\xi:M\to BG$ is a classifying map of $\xi$.

\smallskip

We will also use the projective version of Thom polynomials (see \cite{forms}), as follows.
Assume that $G$ acts on $V$  in such a way that
the scalars are in the image of $G\to GL(V)$. Then the orbits of this action
(different from $\{0\}$) are in bijection with the orbits of the induced
action of $G$ an $\p V$. Also, the corresponding orbits, $\eta$ and
$\p \eta$ have the same codimension. The equivariant Poincar\'e dual
of $\p \eta$ will be called the {\sl projective Thom polynomial} of
$\eta$:
$$\p\tp_\eta=\text{Poincar\'e dual of }\ [\p\eta]\in H^*_G(\p V,\Z)=
H^*(BG,\Z)[x]/(Q(x))  \  \  \  (\deg(x)=2),$$
where $Q(x)$ is the product of all $(x+\alpha_j)$'s, where $\alpha_j\in
H^2(BG)$ are the weights of the representation of $G$ on $V$ \cite{botttu}.
The projective Thom polynomial can be written as
$\p\tp_\eta=p_c+p_{c-1}x+\ldots+p_0x^c$, where $p_i\in H^{2i}(BG)$.
By  \cite[Section 6]{forms},  $p_c=\tp_\eta$ and $p_0$ is the degree of the variety
$\p \eta$.
Seemingly, the projective Thom polynomial contains more information then
the ``affine'' one. This is not the case: $\p\tp_\eta$ can be obtained
from $\tp_\eta$ by a simple substitution, see Theorem~6.1 in \cite{forms}
(although this fact will not be used in the present paper).
In particular,  the degree $p_0$ of $\p \eta$ itself can be obtained from
$\tp_\eta$  by  a  substitution.  For this substitution in our specific case,
see \ref{relations}(2).

\section{Coincident root loci} \label{crl}

Consider the $d$-th symmetric power $V_d=S^d\C^2$ of the standard representation
of $G=GL_2(\C)$,   and the invariant subvariety $X_\lambda$ associated with a  partition
$\lambda=(\lambda_1,\lambda_2,\ldots,\lambda_n)$ of $d$   (cf.  introduction).
In this section we compute  its Thom polynomial  $\tp_{\lambda}\in H^*(BG,\Z)$.

Points in the projectivization $\p V_d$ of $V_d$ can be identified
with $d$-tuples of points in $\p^1=\{(x:y)\}$ (counted with
multiplicities). The projectivization $\p X_{\lambda}$ is then the
closure of the set  of $d$-tuples having  $n$ distinct points with
multiplicities $\lambda_1, \lambda_2,\ldots,\lambda_n$. The
variety $\p X_{\lambda}$ is called the coincident root locus.

Consider also the other notation
$\lambda=(1^{e_1}2^{e_2}\ldots r^{e_r})$  with
$\sum ie_i=d$ and $\sum e_i=n$ (cf. introduction).  Then
$\p X_{\lambda}$ is the image of the map
$$\phi:\p V_{e_1}\times\p V_{e_2}\times \ldots \times \p V_{e_r} \to
\p V_d$$   defined (via point-tuples of $\p^1$) by
$(D_1,D_2,\ldots,D_r)\mapsto \sum iD_i$. It is readily seen that
$\phi$ is birational onto its image $\p X_\lambda$ (i.e. is a resolution
of $\p X_\lambda$).   In particular,  $\dim(\p X_\lambda)=n$ and
$\tp_\lambda$  is  of degree $d-n$ (cf. section \ref{tp}).

The map $\phi$ is equivariant under the action of $G$ on the two
spaces, hence it makes sense to talk about the maps $\phi^*$
(induced by $\phi$) and $\phi_!$ (the push-forward map of $\phi$)
in $G$-equivariant cohomology. The equivariant cohomology rings are
as follows (cf. e.g. \cite{botttu}, p. 270):
$$H^*_G(\,{\textstyle \prod_i}\,
 \p V_{e_i},\Z)=R[x_1,\ldots,x_r]/(Q_{e_1}(x_1),\ldots,
Q_{e_r}(x_r)), \ \mbox{resp.} \ \ H^*_G(\p V_d,\Z)=R[x]/(Q_d(x)).$$
Here
$$R=H^*(BG,\Z)=\Z[c_1,c_2]=\Z[u,v]^{{\mathbb Z}_2},$$
where
${\mathbb Z}_2$ permutes the roots $u$ and $v$ (hence $c_1=u+v$ and
$c_2=uv$);  and the
polynomial $Q_k$ ($k\geq 1$) is  defined by
$$Q_k(y)=\prod_{\alpha \text{\ is a weight of\ } S^k\C^2}
(y+\alpha)=\prod_{j=0}^k(y+ju+(k-j)v).$$
The map $\phi^*$ is a ring homomorphism, it leaves elements of $R$ invariant,
and it maps $x$ to
$$\phi^*(x)=\sum_{i=1}^r ix_i.$$

The above rings can be described also as finite dimensional modules over $R$,
spanned by $\prod_i x_i^{k_i} (0\leq k_i\leq e_i)$ and $x^k (0\leq k\leq d)$,
respectively. A representative of an element  $[f]$  (in any of these rings) is
{\em reduced} if it is written
as an $R$-linear combination of these monomials. It is denoted by $[f]_{red}$.  In this language,
the value of the integration maps (along the fibers)
$$\int_{\prod \p V_{e_i}} : H^*_G(\,{\textstyle \prod_i \p V_{e_i}},\Z) \to
R, \ \ \mbox{resp.} \ \  \int_{\p V_d} : H^*_G(\p V_d,\Z ) \to R$$
are the coefficients of the top degree
monomials in the corresponding reduced forms: i.e. the coefficient of
${\mathbf x}^e:=\prod_i x_i^{e_i}$ in the first case, and the coefficient of
 $x^d$ in the second case.

Set  $q(x):=(Q_d(x)-C_{d+1})/x=x^d+C_1x^{d-1}+\ldots+C_d$, where
$Q_d(x)=\sum_{j=0}^{d+1} C_{d+1-j}(c_1,c_2)x^{j}$ and  $C_0=1$.

\begin{theorem} \label{alg} \
$\tp_\lambda$ equals $\int_{\prod\p V_{e_i}} \phi^*(q)$.
\end{theorem}

\begin{proof}
First we prove that $\tp_\lambda=\int_{\p V_d} (q\cdot \ptp_\lambda)$.
Indeed, from the general theory of projective and affine Thom polynomials
(cf. section \ref{tp}) we know that $\ptp_{\lambda}=p_{d-n}+p_{d-n-1}x+\ldots
+p_0x^{d-n}$, where $p_j\in \Z[c_1, c_2]$ and
$p_{d-n}=\tp_{\lambda}$. When we multiply $x^jp_{d-n-j}$ ($1\leq j\leq d-n$)
with $q=(Q_d-C_{d+1})/x$ and reduce it modulo $Q_d(x)$,
the coefficient of $x^d$ will be 0. So the only contribution
comes from $qp_{d-n}$, which is the coefficient $p_{d-n}$ of $p_{d-n}x^d$.

Now, using the definition of $\ptp_{\lambda}$ and the fact that
$\phi$ is birational, we have
$\tp_{\lambda}=\int_{\p V_d} (q\cdot \phi_!(1))$ which equals
$\int_{\prod\p V_{e_i}} \phi^*(q)$, what we wanted to prove.
\end{proof}

Theorem \ref{alg} gives the following computational recipe: {\em
$Tp_\lambda$ is the top coefficient (i.e. the coefficient of
${\mathbf x}^e$) of  $\phi^*(q)_{red}$.} Notice that any
representative $[f]_{red}$ is automatically computed by computer
algebra packages (e.g. \cite{m2}), hence  one gets an algorithmic
solution of finding the Thom polynomials, see e.g.
www.unc.edu/\~{}rimanyi/progs/rootloci.m2. We can, however, give
explicit formulae as well.

\vspace{3mm}

\noindent {\bf Formulae for Thom polynomials.}

\vspace{3mm}

\begin{lemma}\label{l1} Set $f\in R[y]$ with class $[f]$ $mod\ Q_e(y)$. Then
the top coefficient of the reduced representative $[f]_{red}$ is
\begin{equation*}
\int_{\p V_e}\, [f]=\frac{1}{(v-u)^e}\sum_{s=0}^e \frac{(-1)^{s}
f\big(-(e-s)u-sv\big)}{s! (e-s)!}.\tag{1}\end{equation*}
\end{lemma}

\begin{proof} This is a simple application of the Atiyah-Bott
integration formula \cite[p.9]{ab-moment} but we prefer
to give a direct proof as follows. The formula is linear in $f$,
hence it is enough to verify
it for any $f(y)=y^j$ ($j\geq 0$). In this case
we need $A_e$ where  $y^j\equiv A_ey^e+A_{e-1}y^{e-1}+\ldots
+A_0$ modulo the ideal $(Q_e(y))$.
If we consider this congruence for $y=-eu,-(e-1)u-v,\ldots,-ev$ then we
get a system of equations for $A_e,\ldots, A_0$ (since  $Q_e(y)$
vanishes  at these points). The matrix of this system is
a Vandermonde matrix, so by Cramer's rule we get the formula.
\end{proof}

\begin{corollary}\label{egy}{\bf (The ``naive'' formula)}
 Let $\sum_{s_1,\ldots, s_r}$ denote
the sum over $0\leq s_i\leq e_i$ for each $1\leq i\leq r$. Then

$$\tp_{\lambda}=\frac{1}{(v-u)^n}\sum_{j=n}^d \ \sum_{j_1+\ldots+j_r=j}C_{d-j}
{j \choose j_1,\ldots, j_r} \sum_{s_1,\ldots,s_r}
 \prod_{i=1}^r \frac{(-1)^{s_i} (-i)^{j_i}((e_i-s_i)u+s_iv)^{j_i}}
{s_i!(e_i-s_i)!}.$$
\end{corollary}

\begin{proof} Write $\sum_{j=n}^{d} C_{d-j}(\sum ix_i)^j$ as a linear
combination of monomials of type $\prod_{i=1}^r x_i^{j_i}$.
The polynomial $Q_{e_i}(x_i)$
only contains the variable $x_i$. Hence to find the top coefficient of
the remainder of $\prod_i x_i^{j_i}$, we can simply multiply
the top coefficient of the remainders of $x_i^{j_i}$ modulo $Q_{e_i}(x_i)$.
Therefore, the formula follows from lemma \ref{l1} applied
one-by-one for each $x_i^{j_i}$.
\end{proof}

One can get a more interesting formula as follows. First notice that
$xq+C_{d+1}=Q_d$, hence $(\sum ix_i)\phi^*(q)\equiv
-C_{d+1}$  modulo the ideal ${\mathcal I}\subset R[x_1,\ldots,x_r]$
generated by all $Q_{e_i}(x_i)$ ($1\leq i\leq r$).
We consider the following identities regarding
$1/\sum_i ix_i$. Let $t$ be a free variable. Then
$$\frac{1}{-t+\sum_i ix_i}=
\frac{1}{-t}\sum_{j\geq 0}\  (\,{\textstyle \sum_i} ix_i/t\, )^j=
\frac{1}{-t}\sum_{j\geq 0} \ \sum_{j_1+\ldots+j_r=j}
{j \choose j_1,\ldots, j_r}   \, {\textstyle \prod_i}(ix_i/t)^{j_i}.$$
By lemma \ref{l1}, the top coefficient of the last expression is
$$\frac{1}{-t}\sum_{j\geq 0} \ \sum_{j_1+\ldots+j_r=j}
{j \choose j_1,\ldots, j_r} \, \prod_i
\sum_{s_i=0} ^{e_i} \,  \frac{(-1)^{s_i}}{(v-u)^{e_i}  s_i!(e_i-s_i)!}
\big(\, (e_i-s_i)u+s_i v\,\big)^{j_i}(-i/t)^{j_i}$$
$$=\frac{1}{(-t)(v-u)^n}
\sum_{s_1,\ldots,s_r} \frac{(-1)^{\sum_i s_i}}{\prod_i s_i!(e_i-s_i)!}\cdot
\frac{1}{1+\sum_i i(\, (e_i-s_i)u+s_iv\, )/t}$$
$$=\frac{1}{(v-u)^n}
\sum_{s_1,\ldots,s_r} \frac{(-1)^{\sum_i s_i}}{\prod_i s_i!(e_i-s_i)!}
\, \cdot \, \frac{1}{-t- du + (\sum_i is_i)(u-v)}.$$
Let $A(t)$ be this last expression. The above identities show the following
congruence (valid for generic $t$):
\begin{equation*}
(-t+{\textstyle\sum_i} ix_i)(\, A(t)x^e+\,\mbox{lower order terms}\,)
\equiv 1\ (mod \ {\mathcal I}).
\tag{1}
\end{equation*}
Evidently, this is true for $t=0$ as well. On the other hand, notice that
there is a unique {\em reduced} $Y\in R[x_1,\ldots,x_r]$ satisfying
$(\sum ix_i)Y\equiv -C_{d+1}\ (mod\ {\mathcal I})$. Indeed, if both
$Y$ and $Y'$ satisfy it, then $-C_{d+1}Y'\equiv Y(\sum ix_i)Y'\equiv -
C_{d+1}Y$, hence $Y=Y'$. Since $(\sum ix_i)\phi^*(q)\equiv -C_{d+1}$, from
(1) (with $t=0$) we get that the top coefficient of $\phi^*(q)_{red}$ is
$-C_{d+1}A(0)$. Hence, we proved:

\begin{theorem} \label{ketto} With the notation
$C_{d+1}:=C_{d+1}(S^d\C^2)=\prod_{j=0}^d(\, ju+(d-j)v\,)$, one has
$$\tp_{\lambda}= \frac{C_{d+1}}{(v-u)^n}\cdot
\sum_{s_1,\ldots,s_r} \frac{(-1)^{\sum_i s_i}}{\prod_i s_i!(e_i-s_i)!}
\, \cdot \, \frac{1}{du - (\sum_i is_i)(u-v)}.$$
\end{theorem}
\noindent
This can also be considered as a higher order divided difference formula,
cf. \ref{e2}.

\begin{example}\label{e1} If $\lambda=i^{e_i}$, hence $d=ie_i$, then
$$\tp_\lambda=i^{e_i}\cdot \
\prod_{0\leq j\leq d;i\nmid j}\ (\,ju+(d-j)v\,).$$
\end{example}
\noindent
This can be deduced from \ref{ketto}   (cf.  with the next remark),
but one also can argue  as
follows. Since $ix_i\phi^*(q)+C_{d+1}\equiv 0\ (mod\ Q_{e_i}(x_i))$, clearly
$ix_i\phi^*(q)_{red}+C_{d+1}\equiv 0$ as well. Since $ix_i\phi^*(q)_{red}+
C_{d+1}$ and $Q_{e_i}(x_i)$ both have degree $e_i+1$, one  gets that
$ix_i\phi^*(q)_{red}+C_{d+1}=C\cdot Q_{e_i}(x_i)$ for some $C\in R$.
Comparing the
coefficients of $x_i^{e_i+1}$ and $x_i^0$, one obtains
$$i\tp_\lambda=C_{d+1}(S^d\C^2)/C_{e_i+1}(S^{e_i}\C^2).$$

\begin{remark}\label{r} Lemma \ref{l1} has the following consequence.
For some $C\in R$ and $g\in R[y]$, we denote by $[C/g]_{red}$
(or by $\int_{\p V_e}[C/g]$) that
reduced element which satisfies $[C/g]_{red}\cdot g\equiv C\ (mod\ Q_e(y))$
(if it exists). Then one also has:
\begin{equation*}
\int_{\p V_e}\, [C/g]=\frac{1}{(v-u)^e}\sum_{s=0}^e \frac{(-1)^{s}}
{s! (e-s)!}\cdot \frac{C}{g \big(-(e-s)u-sv\big)}.\tag{1}\end{equation*}
Its proof is similar to the proof of \ref{ketto},  which, in fact, is
a multivariable version of (1)   (applied for $-C_{d+1}/\sum ix_i$).

Let us consider again $\lambda=i^{e_i}$. Theorem \ref{ketto} and (1)
gives that $\tp_\lambda=\int_{\p V_{e_i}} [-C_{d+1}(S^d)/ix_i]$.
But $x_i(x_i^{e_i}+\cdots )+C_{e_i+1}(S^{e_i})=Q_{e_i}$,
hence $\int_{\p V_{e_i}} [-C_{e_i+1}(S^{e_i})/x_i]=1$.
In particular, $\tp_\lambda=C_{d+1}(S^d)/iC_{e_i+1}(S^{e_i})$, as it was
verified in \ref{e1}.
\end{remark}

\begin{example}\label{e2} Assume that $\lambda=i^{e_i}j^{e_j}$ ($i\not=j$).
Consider the expression given by
\ref{ketto} for this $\lambda$, and apply in variable $x_i$
the identity \ref{r}(1). Clearly
$du-(is_i+js_j)(u-v)=g(-e_iu+s_i(u-v))$, where $g(x_i):=a-ix_i$ with
$a:=je_ju-js_j(u-v)$. Therefore
$$\tp_\lambda=
\frac{C_{d+1}(S^d)}{(v-u)^{e_j}}\sum_{s_j=0}^{e_j} \frac{(-1)^{s_j}}
{s_j! (e_j-s_j)!}\cdot \int_{\p V_{e_i}}\, [1/g(x_i)].$$
Since $Q_{e_i}(x_i)-Q_{e_i}(a/i)=(x_i-a/i)(x_i^{e_i}+\cdots)$ one gets
$\int_{\p V_{e_i}} [iQ_{e_i}(a/i)/g(x_i)]=1$. Hence
\begin{equation*}\tp_\lambda=
\frac{C_{d+1}(S^d)}{(v-u)^{e_j}}\sum_{s_j=0}^{e_j} \frac{(-1)^{s_j}}
{s_j! (e_j-s_j)!}\cdot \frac{1}{i\cdot Q_{e_i}(\, (je_ju-js_j(u-v)\,)/i)}.
\end{equation*}
%
For example, assume that $\lambda=i^{e_i}j$, i.e. $e_j=1$. Then
$s_j=0$ or $1$, hence
$$\tp_\lambda=
\frac{C_{d+1}(S^d)}{i(v-u)}\cdot \Big(
\frac{1}{Q_{e_i}(ju/i)}-\frac{1}{Q_{e_i}(jv/i)}\Big).$$
It is convenient to express this in the language of {\em divided difference}:
If $P(u,v)$ is a polynomial in two variables $(u,v)$, we denote by
$\dd(P)$ the polynomial $(P(u,v)-P(v,u))/(u-v)$. Then
$$\tp_{(i^{e_i}j)}=\frac{1}{i}\cdot \dd\, \Big(\, \frac{C_{d+1}(S^d)}
{Q_{e_i}(jv/i)}\, \Big)\, =
i^{e_i}\cdot \dd\, \Big(\prod \big(\, (d-k)v+ku\, \big)\Big),$$
where the product is over $k$ with $0\leq k\leq d$, but $k\not=is$ with
$0\leq s\leq e_i$.
In particular,
\begin{equation*}
\tp_{(1^{e_1}j)}=
\dd\, \Big(\prod_{l=0}^{j-1} \big(\, lv+(e_1+j-l)u\, \big)\Big) \ \ (\mbox{for} \ j\geq 2),
\tag{1}
\end{equation*}
which is  equivalent with  Kirwan's formula  \cite[page 902]{kirwancikk}.
\end{example}

\begin{example}\label{e3} Assume that $d=2h$ is even, $h>2$ and $\lambda=
(1^{h-j},j,h)$ for some $1<j<h$.
By a similar argument as in \ref{e2} and by a computation,
one has
$$\tp_\lambda=
\frac{C_{d+1}(S^d)}{(u-v)^2}\cdot \Big(
\frac{1}{Q_{h-j}(hu+ju)}-
\frac{1}{Q_{h-j}(hu+jv)}-
\frac{1}{Q_{h-j}(hv+ju)}+
\frac{1}{Q_{h-j}(hv+jv)}
\Big)$$
$$=\dd\Big[\ \frac{C_{d+1}(S^d)}{u-v}\cdot \Big(
\frac{1}{Q_{h-j}(hv+jv)}-
\frac{1}{Q_{h-j}(hv+ju)}\Big) \Big]
=\dd\,\Big[ D_j\cdot \prod_{l=0}^{h-1}\big(\,lv+(d-l)u\,\big)\Big],$$
where
$$D_j:=\frac{1}{u-v}\cdot\Big[
\prod_{l=h-j+1}^{h}(\,lv+(d-l)u\,)- \prod_{l=0}^{j-1}(\,lv+(d-l)u\,)
\Big].$$

E.g., if $j=2$, then
\begin{equation*}
\tp_{\lambda}=h(h-1)\cdot
\dd\,\Big[ (u+3v)\prod_{l=0}^{h-1}\big(\,lv+(d-l)u\,\big)\Big].
\end{equation*}
\end{example}
\begin{remarks}\label{relations} \ {\em (1)
The Thom polynomials are connected by many interesting polynomial  relations.
E.g., the next section presents two situations when
the ideal generated by natural families of Thom polynomials is
generated only by two of them.  Some of these relations can be verified easily.
E.g., assume $d=2h$ as in \ref{e3}, consider the
partitions $\lambda_0'=(1^{h-2},2,h)$,
$\lambda_0=(1^h,h)$, $ \lambda_1=(1^{h-1},h+1)$ and
$\lambda_2=(1^{h-2},h+2)$. Then from
\ref{e2}(1) and \ref{e3}, one  gets
$\tp_{\lambda_1}=hc_1\cdot \tp_{\lambda_0}$ and
$$(h-1)\cdot \tp_{\lambda_2}=
(h-1)(h-2)c_1\cdot  \tp_{\lambda_1}+c_1\tp_{\lambda_0'}.$$
(2)  Using \cite{forms}, one can determine  $\deg(\p X_\lambda)$ by the substitution
$u=v=1/d$ in $\tp_\lambda\in \Z[u,v]$.  The interested reader is invited to verify
the compatibility of Hilbert's result (cf. introduction)  with this section.

\noindent (3)  In the sequel we will use many times the following divided difference formula.
For any polynomial $A\in \Q[u,v]$ write $A^*(u,v):=A(v,u)$. Then }
\begin{equation*}
\dd(AB)=B^*\cdot \dd (A)+A\cdot \dd (B).\end{equation*}
\end{remarks}

\section{Thom polynomial description of the cohomology ring of the moduli space}\label{moduli}
In this section we apply the coincident root loci formulas in the study of the cohomology ring of the moduli space of the representation $S^d\C^2$ (in the Geometric Invariant Theory sense). We calculate the rational cohomology rings $H^*_G(X^{ss})$, $H^*(X^{ss}/G)$ and  $H^*_G(X^{s})\cong H^*(X^{s}/G)$ in terms of generators and relations. If $d$ is odd then all these rings coincide, but for the even case they are different.

There is an extensive literature on these cohomology rings,
both from combinatorial-algebraic (see e.g. \cite{brion},
\cite{martin}) and from geometric point of view
(the Atiyah-Bott-Kirwan theory \cite{kirwan}).
Our approach (in the odd $d$ case) is closest to that of Kirwan.
The advantage of our approach is that we treat the odd and even cases
in a uniform language, and that we provide for the above
cohomology rings a very transparent structure: we obtain
explicit presentations of them in terms of
generators and relations with clear geometric meanings.

Let us sketch our approach in the odd case first (for details see below). In this case the Kirwan stratification of $S^d\C^2$ is $G$-perfect since the normal (equivariant) Euler classes of the strata are not zero-divisors. It implies that the spectral sequence of the corresponding filtration degenerates. It is not difficult to calculate all but
the 0\textsuperscript{th} column of the $E_1$-table, so by subtraction we can calculate the ranks of the  0\textsuperscript{th} column: the Betti numbers of $H^*_G(X^{ss})$. Also by $G$-perfectness the natural map
$$\kappa:H^*_G(S^d\C^2)\cong \Q[c_1,c_2] \to H^*_G(X^{ss})$$
is surjective so we need to find relations in terms of $c_1$ and $c_2$, ie. we have to find generators of $\Ker (\kappa)$. If $Y\cap X^{ss}=\emptyset$ for an invariant subvariety $Y$ then clearly $[Y]\in \Ker (\kappa)$. (This idea was studied in \cite{ss}). So all the higher Kirwan strata provide relations. But the Kirwan strata are coincident root loci for specific partitions and we can calculate their equivariant Poincar\'e dual using the first part of the paper. It turns out that the first two Kirwan strata are enough to generate  $\Ker (\kappa)$ which can be checked by a simple Betti number calculation.

The main difficulty in the even case is that for one of the strata in a refined Kirwan stratification the normal (equivariant) Euler class is a zero-divisor. To prove $G$-perfectness we use the results of the first part of the paper. Namely we show that certain elements in the $E_1$-table can be represented by the Poincar\'e dual of coincident root loci (these are not  Kirwan strata!) and they survive to $E_\infty$, hence they could not be hit by a differential. After $G$-perfectness is proven the process is the same as in the odd case. We can find  coincident root loci in the null cone such that their Poincar\'e dual generate  $\Ker (\kappa)$. Here we also need two coincident root loci but one of them is not  a Kirwan stratum.

\smallskip

In this section all cohomologies are meant with rational coefficients.

Let us consider the  Kirwan-stratification (see \cite{kirwancikk} and \cite{kirwan}) of the vector space $V_d$:
\begin{itemize}
\item
$X^{ss}=\{B\ |\ B\ \text{has no root of multiplicity} >d/2 \},$
\item
$X_{i}=\{B \ |\ B\ \text{has a root of multiplicity}\ i\ \text{but no with multiplicity}\ i+1\}$ ($d/2<i\leq d$),
\item
$X_0=\{0\}$.
\end{itemize}
The strata are smooth open submanifolds, the complex codimensions are
$0,\ i-1,\ d+1$ in the three cases.   By $F_i=\cup$ strata of complex
codimension $\leq i$ we get a filtration of $V_d$:
$$\emptyset=F_{-1}\subset F_0 \subset F_1 \subset \ldots \subset
 F_{d+1}=S^d\C^2.$$
Let $E_*^{*,*}$ be the associated spectral sequence in $G$-equivariant
cohomology with $\Q$-coefficients.

\begin{proposition}\label{specseq}\
\begin{enumerate}
\item $E_1^{0,*}=H_G^*(X^{ss};\Q)$
\item $E_1^{2p,*}=H^*(BU(1);\Q)$ for $p=[d/2],\ldots,d-1$;
\item $E_1^{2(d+1),*}=H^*(BG;\Q)$;
\item $E_1^{*,*}=0$ for all cases not covered by (1), (2), (3);
\item The spectral sequence converges to $H^*(BG;\Q)$;
\item The spectral sequence degenerates at $E_1^{*,*}$ (in particular,
 $H^{odd}_G(X^{ss},\Q)=0$).
\end{enumerate}
\end{proposition}

\begin{proof} By definition we have $E^{2p,*}_1=H_G^{2p+*}(F_p,F_{p-1})$
which is by
Thom isomorphism $H^*_G(F_p\setminus F_{p-1})$. This proves (1) and (4).
For $p=d+1$ we have $E_1^{2(d+1),*}=H^*_G(\{0\})=H^*(BG)$ which proves (3).
For $d/2<i\leq d$ we define $Y_i=\{B\in X_i \, : \, x^i|B \ \mbox{and} \
${\em coeff}$\, (x^iy^{d-i})=1\}$.
 Let $H$ be the stabilizer subgroup of
$Y_i$, i.e. the group of matrices of the form
$\begin{pmatrix} \alpha_1 & \beta \\ 0 & \alpha_2 \end{pmatrix}$
with $\alpha_1^i\alpha_2^{d-i}=1$.
Since $Y_i$ is contractible, and  $X_i=G\times_H Y_i$, part (2) follows from
$$H^*_G(X_i)\cong H^*_G(G\times_H Y_i)\cong H^*_H(Y_i)\cong H^*(BH)\cong
H^*(BU(1)) \ (\mbox{over $\Q$}).$$
The degeneracy of the spectral sequence---called $G$-perfectness by
Atiyah-Bott in \cite{ab}---follows from usual arguments, as follows.
Let us build up $V_d$ by gluing the strata one by one together in order
of increasing codimension. Then at one step we have $U$ and glue a new
stratum $X$ of complex codimension $c$ to it.
We need to prove that the first map in the diagram
$$H^{n-2c}_G(X)\cong H^n_G(U\cup X,U)\to H^n_G(U\cup X) \to  H^n_G(X)$$
is injective. However, the whole composition is the multiplication with the
equivariant Euler class of the stratum $X$. This is an injective map
being a multiplication by a non-zero element in a polynomial ring.
(For a computation of an equivariant Euler class see the proof of \ref{pr2}.)
\end{proof}

Since $E_{\infty}=E_1$, the sum of the ranks of the groups in diagonal
(i.e. $p+q=r$) entries must be the rank of the appropriate cohomology group
of $H^*(BG;\Q)$. Thus we have the following

\begin{corollary}\label{hs} Let $h:=[d/2]$.
The Poincar\'e  series of the ring $H^*_G(X^{ss};\Q)$ is
$$\frac{1}{(1-t)(1-t^2)}(1-t^{d+1})-\frac{1}{1-t}(t^{h}+\ldots+t^{d-1})=
\frac{1-t^h-t^{h+1}+t^d}{(1-t)(1-t^2)} \ \ (\deg(t)=2).$$
\end{corollary} \qed

What we obtained so far is basically equivalent to the Atiyah-Bott-Kirwan
theory applied to our representation, see \cite[16.2]{kirwan}.

What can also be seen from the spectral sequence is that
$H^*_G(X^{ss})=H^*(BG)/I$ where the ideal comes from
the $p>0$ columns of the spectral sequence. Thus among the
elements of $I$ we have the ones that are the images of the
generators of $E_1^{2p,0}$ under the edge-homomorphism.
For $[d/2]\leq p\leq d-1$, these are exactly
the Thom polynomials corresponding to the strata $X_i$, $i=p+1$.
We have $\tp(X_i)=\tp_\lambda$ whit $\lambda=(1^{d-i},i)$, since
the  closures  of $X_i$ and $X_\lambda$ are the same. The above Betti number computation
can be used to test if a few of these Thom polynomials are enough
to generate $I$.

\begin{theorem}\label{theo} Set
 $\lambda_1=(1^{d-h-1},h+1)$ and  $\lambda_2=(1^{d-h-2},h+2)$, where $h=[d/2]$. Then
$I$ is generated by $\tp_{\lambda_1}$ and $\tp_{\lambda_2}$. In particular,
$$H_G^*(X^{ss};\Q)=\Q[c_1,c_2]\Big/ (\tp_{\lambda_1},\tp_{\lambda_2}).$$
\end{theorem}

\begin{proof} We already  observed that the given two
$\tp$'s are in $I$. Now we prove that the ring on the right hand side
has the same Poincar\'e series as the one given in Corollary \ref{hs}.

We claim  that  the ideal $J:=(\tp_{\lambda_1},\tp_{\lambda_2})$
 has the following $R$--resolution:
$0\leftarrow J \leftarrow R(h)\oplus R(h+1) \leftarrow R(d)\leftarrow 0$.
If $d=2h+1$ then for this we only need to prove
that $\tp_{\lambda_1}$ and $\tp_{\lambda_2}$
have no nontrivial common divisor $D$. We know that
$\tp_{\lambda_1}=\dd(\Pi)$, $\tp_{\lambda_2}=\dd(\Pi L)$, where
$\Pi(u,v)=\prod_{l=0}^{h} (lv+(d-l)u)$ and $L(u,v)=(h+1)v+hu$.
By \ref{relations}(3), if $D|gcd(\tp_{\lambda_1},\tp_{\lambda_2})$, then
$D|\Pi$, hence $D|gcd(\Pi,\dd(\Pi))$ as well. But $gcd(\Pi,\Pi^*)=1$,  which ends the
proof of the claim.

So we get the Poincar\'e  series of $R/J$ as
$(1-t^h-t^{h+1}+t^{2h+1})/(1-t)(1-t^2)$, which is the same as the
Poincar\'e  series of $H_G^*(X^{ss};\Q)$. For $d$ even the proof is similar.
\end{proof}

\begin{discussion}\label{ss} {\bf The cohomology ring  of $X^{ss}/G$.}
Observe that if $d$ is odd then $X^{ss}=X^s$, and
all stabilizers of polynomials in $X^{ss}$ are
finite. Therefore, we have
the ring isomorphism $H_G^*(X^{ss};\Q)=H^*(X^{ss}/G;\Q)$
with Poincar\'e polynomial $(1-t^h)(1-t^{h+1})/(1-t)(1-t^2)$.

If $d=2h$ is even, then $X^{ss}/G=X^s/G\cup\{p^{ss}\}$, where $p^{ss}$ is the
unique ``semisimple point'' of $X^{ss}/G$.
The Poincar\'e series of $H^*_G(X^{ss})$ is infinite; it is:
\begin{equation*}
\frac{1}{1-t^2}+t\cdot P(t), \ \mbox{where $P(t)$ is the polynomial}
\ \frac{(1-t^{h-1})(1-t^h)}{(1-t)(1-t^2)}\ \ \ (\deg(t)=2).
\tag{1}
\end{equation*}
All the stabilizers of the
stable part are finite, and there is only one orbit in the strict
semistable part with infinite stabilizer $H^{ss}$,
namely the orbit of the partition $(h,h)$. $H^{ss}$ can be described
explicitly, and one has an exact  sequence $1\to U(1)\times \Z_h
\to H^{ss}\to \Z_2\to 1$. Hence $BH^{ss}$ is a double
covering of $BU(1)\times B\Z_h$ with rational cohomology $H^*(BH^{ss})=
H^*(BU(1))^{\Z_2}=
\Q[t]^{\Z_2}$ ($\deg_t=2$). Here the $\Z_2$-action is $t\mapsto \pm t$, hence
the invariant part is $\Q[t^2]$ with an infinite Poincar\'e series
$1/(1-t^2)$. This is exactly the ``infinite contribution'' in the above
Poincar\'e series of $H^*_G(X^{ss})$.

In fact, the map $r:H^*(BG)\to H^*(BH^{ss})$ (induced by the inclusion)
is the following.  At the level of roots, it is given by
$u\mapsto \pm t$ and $v\mapsto \mp t$, hence it is the epimorphism
$r:\Q[c_1,c_2]\to \Q[t^2]$ given by $c_1\mapsto 0$ and $c_2\mapsto -t^2$.

As usual, for any connected space $Z$,  let $\tilde{H}^*(Z)$ be
the kernel of $H^*(Z)\to H^*(\mbox{point})$, as an ideal (or subring without
unit) in $H^*(Z)$. The ring  $H^*(Z)$ can be
reconstructed from $\tilde{H}^*(Z)$ by adding the unit: $H^*(Z)=\Q\langle 1
\rangle \oplus \tilde{H}^*(Z)$ (with the natural multiplication).

Let $o$ be the orbit corresponding to the partition
 $(h,h)$ and consider the natural inclusion
$j:o\times_GEG \to X^{ss}\times_GEG$. Obviously, $o\times_GEG$
can be identified with $BH^{ss}$.
Moreover, $j^*:H^*_G(X^{ss})\to H^*(BH^{ss})$ induced  by $j$
can be identified with the epimorphism
$\Q[c_1,c_2]/(\tp_{\lambda_1},\tp_{\lambda_2})\to
\Q[t^2]$, $c_1\mapsto 0$ and $c_2\mapsto -t^2$ induced by $r$ above.
In fact, $\tp_{\lambda_1}$ and $\tp_{\lambda_2}$ are both divisible by
$c_1$ (cf. \ref{relations}(1)), hence $r$ sends the ideal generated by them to zero.

Finally, notice that $H^*(X^{ss}\times_GEG,BH^{ss})=
\tilde{H}^*(X^{ss}\times_GEG/
BH^{ss})$, and the natural map
$r:X^{ss}\times_GEG/BH^{ss}\to X^{ss}/G$ induces an
isomorphism at the level of rational cohomology rings. In particular,
the long exact cohomology sequence
of the pair $(X^{ss}\times_GEG,BH^{ss})$ transforms into the short
exact sequences:
\begin{equation*}
0\to \tilde{H}^*(X^{ss}/G)\to H^*_G(X^{ss})\stackrel{j^*}{\to} H^*(BH^{ss})
\to 0.\tag{2}
\end{equation*}
Analyzing the kernel of $j^*$, we get:
\end{discussion}

\begin{corollary}\label{corss} With the notations of \ref{theo}, one
has the following ring isomorphisms:
$$H^*(X^{ss}/G;\Q)=\Q[c_1,c_2]\Big/ (\tp_{\lambda_1},\tp_{\lambda_2})
\ \ \mbox{if $d$ is odd;}$$
$$H^*(X^{ss}/G;\Q)=\Q\langle 1\rangle\oplus (c_1\Q[c_1,c_2])\Big/
(\tp_{\lambda_1},\tp_{\lambda_2})
\ \ \mbox{if $d$ is even}.$$
\end{corollary}
\noindent
Notice that the Poincar\'e series formula \ref{ss}(1)
is compatible  with the \ref{ss}(2) and \ref{corss}.
In particular, if $d=2h$, the Poincar\'e polynomial of $H^*(X^{ss}/G)$ is
$1+tP(t)$.

\begin{discussion}\label{s} {\bf The cohomology ring of $X^s/G$.}
Next, for the case $d=2h$,
we wish to determine the cohomology ring of the geometric
quotient $X^s/G$. In the notations below it is convenient to assume $h> 2$
(if $h=2$, then $X^{ss}/G=\p^1$, and $X^s/G=\C$).

We consider a similar spectral sequence, but now  associated with the stratification

\noindent $\bullet$
$X^{s}=\{B\ |\ B \ \text{has no root of multipicity} \geq h \},$

\noindent $\bullet$
$X_{i}=\{B\ |\ B\ \text{has exactly one root of multiplicity}\ i\
\text{but no roots of multiplicity}\ i+1\}$ ($h\leq i\leq d$),

\noindent $\bullet$
$o=\{\text{the orbit associated with the partition $(h,h)$}\}$,

\noindent $\bullet$
$X_0=\{0\}$.

\vspace{2mm}

In lemma \ref{specseq},  $E_1^{0,*}$ will be replaced by $H^*_G(X^s)$.
For $i>h$, the stratum $X_i$ is the same as in the previous case.
But there are two new strata, namely $X_h$ and $o$.
Since $o$ is an orbit with stabilizer $H^{ss}$, $H^*_G(o)=H^*(BH^{ss})$.
The complex codimension of $o$ in $V_d$ is $d-2$, hence this will provide an
additional direct sum contribution in   $E_1^{2(d-2),*}$. Hence, $E_1^{2p,*}=
H^*(BU(1))$ if $h\leq p\leq d-1$, but $p\not=d-2$; and $E_1^{2(d-2),*}=
H^*(BU(1))\oplus H^*(BH^{ss})$. Finally, we compute $E_1^{2(h-1),*}=
H^*_G(X_h)$. Set
$$Y_h=\{B\in X_h \, : \, B=x^h\cdot B'=
x^h(y^h+ a_2x^2 y^{h-2}+\cdots +a_{h}x^h);\
\text{and $B'$ is not an $h$-power}\}.$$

The stabilizer subgroup $H$ of $Y_h$ is the group of diagonal matrices of the
 form $diag(\alpha_1,\alpha_2)$ with $\alpha_1^h\alpha_2^{h}=1$.
One can verify that $X_h=G\times_H Y_h$. Moreover, $B'$ is not a $h$-power
if and only if $(a_2,\ldots,a_h)\not= (0,\ldots ,0)$. Hence $Y_h$ is
$\C^{h-1}\setminus \{0\}$ and the action of $H$ is a diagonal torus
action (modulo a finite group). In particular, $E_1^{2(h-1),*}=
H^*_G(X_h)$ equals the cohomology ring of a weighted projective space of dimension $h-2$,
which is $\Q[t]/(t^{h-1})$ ($\deg(t)=2$).

\begin{proposition}\label{pr2} The spectral sequence converges to $H^*(BG;\Q)$ and it degenerates at $E_1^{*,*}$.
\end{proposition}
\begin{proof} The Euler classes of the strata are not zero-divisors except for $X_h$. So we need the following local version of the Atiyah-Bott argument:
  \begin{lemma}
    Suppose that $\{X_i\}$ is a $G$-equivariant stratification of $V$ and the equivariant normal Euler class of $X_i$ is not a zero-divisor if $\codim(X_i)>c$. Then all differentials of the corresponding spectral sequence $E_r^{p,q}$ starting or landing in the region $p>c$ are zero.
  \end{lemma}
  \begin{proof}[Proof of Lemma] Let $X$ be the union of $X_i$ with $\codim(X_i)>c$. Then the Lemma is equivalent with the statement that $H^*_G(V,V\setminus X)\to H^*_G(V)$ is injective, since $\{E_r^{p,q}:p>c\}$ converges to $H^*_G(V,V\setminus X)$. Injectivity can be proved by adding the $X_i$'s one by one, and noticing that the composition
$$H^{n-2c}_G(X)\cong H^n_G(U\cup X_i,U)\to H^n_G(U\cup X_i) \to  H^n_G(X_i)$$
is   multiplication with the equivariant normal Euler class of the stratum $X_i$ (where $U$ is an open subset of $V$ in which $X_i$ is closed).
  \end{proof}

For the convenience of the reader we show
how one determines the equivariant Euler class of
$o$. Fix an element,
say $x^hy^h$ on $o$, let $H^{ss}$ be its stabilizer, consider an
$H^{ss}$ invariant  normal slice $N$ at $x^hy^h$.
In fact, for $N$ one can take the vector space
spanned by $x^iy^{d-i}$, where $0\leq i\leq d$, but $i\not\in \{h-1,h,h+1\}
$. $H^{ss} $ acts on $N$, and our goal is the computation of the Euler
class $e^{ss}\in H^*(BH^{ss})$
of $EH^{ss}\times_{ H^{ss}}N\to BH^{ss}$. Consider now
the subgroup $U(1)$ of $H^{ss}$ (see \ref{ss}). The Euler class $e\in
H^*(BU(1))=\Q[t]$ of $EH^{ss}\times _{U(1)}N\to BU(1)$ can be computed
as follows. The eigenvalues of $diag(\alpha,\bar{\alpha})\in U(1)$
on $N$ are $(\alpha^d,\alpha^{d-2},\ldots,\alpha^4,\alpha^{-4},\ldots,
\alpha^{-d})$, hence $e=(dt)((d-2)t)\cdots (4t)(-4t)\cdots (-dt)=
mt^{d-2}$ for some $m\not=0$. Since $d$ is even, this is in the
invariants part $H^*(BH^{ss})=\Q[t^2]$ and can be identified in this ring
by $e^{ss}$. Hence $e^{ss}\not=0$.

This type of argument is not working for the stratum $X_h$ (since the
stabilizer of its points are finite, and also $H^*_G(X_h)$
has zero divisors).

In order to show that the differentials $d_{2h-2}^{0,q}$  ($q$ odd and
$2h-3\leq q\leq 4h-7$) of the spectral sequence are trivial,
we consider
another spectral sequence associated with only two strata, namely with
$X^s$ and $X_h$. The differential $d_{2h-2}^{0,q}$ in the two spectral sequences
coincides. If we compare them by the natural maps, then we get the exact
sequence
$$0\to I'\to H^*_G(V_d)\stackrel{\tau}{\to} H^*_G(X^s\cup X_h)$$
where the ideal $I'$ is generated by all the columns  $E_1^{>2h-2,*}$.
In $E_{\infty}^{2(h-1),2(j-1)}$ we can find special
elements, those represented by the Thom polynomials $\tp_j\in H^*_G(V_d)$
associated with the partitions $(1^{h-j},j,h)$, where $0<j<h$. Hence,
$d_{2h-2}^{0,2j+2h-5}=0$ if $\tau(\tp_j)\not=0$, or
equivalently, if $\tp_j\not\in I'$. Notice that the graded ideal $I'$
and the graded ideal $I$ considered in \ref{hs} and \ref{theo} are the
same in the relevant degrees, hence it is enough to verify that $\tp_j \not\in I$ for any $j$. But in \ref{theo} we verified that $I=(\tp_{\lambda_1},\tp_{\lambda_2})$. Hence, we
need to prove:
\begin{equation*}
\tp_j\not\in ( \tp_{\lambda_1},\tp_{\lambda_2}).\tag{1}\end{equation*}
Set
$$\Pi:=\prod_{l=0}^{h-1}(\, lv+(d-l)u\,) \ \ \mbox{and} \ \ L=(h+1)v+(h-1)u.$$
From \ref{e2}(1) one gets $\tp_{\lambda_1}=hc_1\cdot \dd(\Pi)$ and
$\tp_{\lambda_2}=hc_1\cdot \dd(\Pi L)$.
In particular, by \ref{relations}(3),  $\tp_{\lambda_2}=hL^*c_1\cdot \dd(\Pi)-2hc_1\Pi$, hence
\begin{equation*}
( \tp_{\lambda_1},\tp_{\lambda_2})=(c_1\cdot \dd(\Pi), c_1\cdot \Pi).
\tag{2}\end{equation*}
Assume that (1) is not true and we have $\tp_j=Ac_1\cdot \dd(\Pi)+
Bc_1\Pi$. Since the degrees of $\tp_j$  and $\Pi$ are  $h+j-2$ and
$h$ respectively, the degree of $Ac_1$ is $j-1$.
From \ref{e3} and \ref{relations}(3), $\tp_j=\dd(\Pi\cdot D_j)=
D_j^*\cdot \dd(\Pi)+\Pi\cdot \dd(D_j)$. This means that
\begin{equation*} \Pi(\dd (D_j)-Bc_1)=\dd(\Pi)(Ac_1-D_j^*).\tag{3}\end{equation*}
But it is easy to verify that $gcd(\Pi,\dd(\Pi))=1$. Indeed, if
$F|gcd(\Pi,\dd(\Pi))$, then also $F|(u-v)\dd(\Pi)=\Pi-\Pi^*$,
hence $F|\Pi^*$ as well. But $gcd(\Pi,\Pi^*)=1$.

This fact together with (3) show that $\Pi|Ac_1-D_j^*$, but $\deg(Ac_1-D_j^*)
=j-1<\deg\Pi$, hence $Ac_1=D_j^*$. In particular, $c_1|D^*_j$, or $u+v|D_j$.
But this leads to a contradiction. Indeed, analyzing in \ref{e3} the
expression of $(u-v)D_j$, one sees that  the first product
is divisible by $u+v$ (take $l=h$) but the second is not.
Hence, (1) is true.
\end{proof}

\noindent
By similar argument as in the case of $H^*_G(X^{ss})$,  for
$H^*_G(X^s)=H^*(X^s/G)$ one gets:

\begin{corollary}\label{cors}  $H^{odd}(X^{s}/G,\Q)=0$, and the
Poincar\'e series of $H^*(X^s/G)$ is  the polynomial $P(t)$ introduced in
\ref{ss}(1).
\end{corollary}

Let $I''$ be the ideal in $H^*(BG)=\Q[c_1,c_2]$
generated by the columns $E_1^{>0,*}$.
Then one has the ring isomorphism  $H^*(X^s/G)=
\Q[c_1,c_2]/I''$. Now we will consider two special elements of $I''$,
namely the  Thom polynomials
$\tp_{\lambda_0}$ and $\tp_{\lambda_0'}$, where
$\lambda_0=(1^h,h)$ and $\lambda_0'=(1^{h-2},2,h)$. Their degrees
are $h-1$ an $h$ respectively. We will verify now that they are relative
prime. Indeed,  using the above notations,
 $\tp_{\lambda_0}=\dd(\Pi)$ (from \ref{e2}). Moreover, by \ref{e3}
and \ref{relations}(3) one has
$\tp_{\lambda_0'}=h(h-1)\dd((u+3v)\Pi)=h(h-1)[(v+3u)\dd\Pi-2\Pi]$.
In particular, $gcd(\tp_{\lambda_0},\tp_{\lambda_0'})=gcd(\Pi,\dd\Pi)$
which is 1 by the proof of \ref{pr2}. Then the usual Poincar\'e polynomial
argument shows $I''=(\tp_{\lambda_0}, \tp_{\lambda_0'})=(\dd\Pi,\Pi)$.

This can be compared with (2) from the proof of \ref{pr2}:
$(\tp_{\lambda_1},\tp_{\lambda_2})=(c_1\dd\Pi,c_1\Pi)$ (fact which can be
deduced from \ref{relations}(1) as well). Hence we proved:

\begin{theorem}\label{theo2} Assume $d=2h$ and set $\Pi:=\prod_{l=0}^{h-1}(lv+(d-l)u)$. Then:
$$H^*(X^s/G;\Q)=
\Q[c_1,c_2]/(\tp_{\lambda_0},\tp_{\lambda_0'})=
\Q[c_1,c_2]/(\Pi,\dd\Pi),$$
$$H^*(X^{ss}/G,\Q)=\Q\langle 1\rangle \oplus \frac{c_1\Q[c_1,c_2]}{(c_1\Pi,
c_1\dd\Pi)},$$
and the restriction map $H^*(X^{ss}/G)\to H^*(X^s/G)$ is induced by the
identity of $\Q[c_1,c_2]$.
\end{theorem}
\end{discussion}

\begin{discussion}\label{link} {\bf The cohomology ring of the link.}
Denote by $L^{ss}$  the link of the unique semisimple point
$p^{ss}$ in $X^{ss}/G$ (i.e.  $L^{ss}=\rho^{-1}(\epsilon)$,
where $\rho: X^{ss}/G\to [0,\infty)$ is a real analytic map with
$\rho^{-1}(0)=\{p^{ss}\}$ and $\epsilon$ is sufficiently small).
Write   $CL^{ss}$ for the real cone over it (i.e. $CL^{ss}=[0,1]\times L^{ss}/\{0\}\times L^{ss}$).
Then $H^*(CL^{ss},L^{ss})=
H^*(X^{ss}/G, X^s/G)$. Hence $H^*(L^{ss})$ is completely determined by
the restriction morphism from \ref{theo2}.
In fact, $L^{ss}$ is a rational homological manifold of real dimension
$4h-7$ (with Poincar\'e duality). [This can also be proved as follows:
The geometric quotient of the set of ordered $d$-points  of $\p^1$ is
smooth, and one has only finitely many ordered semisimple points.
Hence, $L^{ss}$ is the quotient by a finite permutation group of a smooth
$(4h-7)$--dimensional link].  \ref{theo2},  this duality and a computation  give:

\begin{theorem}\label{theo3} $H^*(L^{ss},\Q)$ can be  generated by two
elements, $c_2$ of degree 4 and $g$ (the Poincar\'e dual of $c_2^{[h/2]-1}$)
of degree $4h-4[h/2]-3$ with relations
$c_2^{[h/2]}=0$ and $g^2=0$. (Notice that all the Betti numbers are 0 or 1.)
\end{theorem}
\end{discussion}

\begin{remark} \ref{theo2} implies the following: the cohomology ring
of the  quasi-projective variety $X^s/G$ of (complex) dimension
$d-3$ shares the Poincar\'e duality properties of a smooth projective
variety of dimension $d-4$. In fact, cohomologically (over $\Q$),
$X^s/G$ behaves like a line bundle ${\mathcal L}$ with Chern class $c_1$
over a smooth projective variety $M$ with cohomology $\Q[c_1,c_2]/(\dd\Pi,
\Pi)$; and $X^{ss}/G$ behaves
like the Thom space of this line bundle (or equivalently,
the complex cone over $M$ associated with ${\mathcal L}$).
In particular, $L^{ss}$ has the cohomology of the $S^1$-bundle of
${\mathcal L}$.
\end{remark}

\begin{remark} Assume that $d=2h+1$ is odd.
It is tempting to compare the moduli space $X^{s}/G$ with
the (possibly weighted) Grassmanian $Gr_2\C^{h+1}$ because
the presentation of their cohomology rings have the same structure
$\Q[c_1,c_2]/ (\dd p_1 , \dd p_2 )$ (where $\deg p_1=h+1$ and
$\deg p_2=h+2$), and they share the  same Betti numbers.
Indeed for the Grassmanian we can take $p_1=u^{h+1}, p_2=u^{h+2}$.
In fact, this analogy can be continued: in both cases
the set of relations are guided by some nice generating function, as follows.
Set $\Pi_0:=1$ and $\Pi_j:=\prod_{l=0}^{j-1}(lv+(d-l)u)$, and consider
the generating function
$${\mathcal G}(q)=\sum_{j\geq 0}{\mathcal G}_jq^j:=
\sum_{j\geq 0}\Pi_jq^j/j!=[1+(u-v)q]^{du/(u-v)}\
\in \Q[u,v][[q]].$$
Then $H^*(X^s/G)=\Q[c_1,c_2]/I$, where $I$ is generated by
$\dd{\mathcal G}_j$, $j>h$.

In the Grassmanian case the same fact is true with
 ${\mathcal G}(q)=1+uq+u^2q^2+\cdots=1/(1-uq)$.\\
However, easy computation shows that, as graded rings, these
cohomology rings are {\em not} isomorphic (unless for small $d$'s).
\end{remark}

\bibliography{root}
\bibliographystyle{alpha}

\end{document}